\documentclass[12pt]{amsart}
\usepackage{epsfig,amsmath,amssymb,amsthm}
\newtheorem{Theorem}{Theorem}[section]

\newtheorem{Corollary}[Theorem]{Corollary}
\newtheorem{Proposition}[Theorem]{Proposition}

\theoremstyle{definition}
\newtheorem{Example}[Theorem]{Example}
\newtheorem{Remark}[Theorem]{Remark}

\begin{document}


\title[Trace Homomorphism]
{Trace Homomorphism for Smooth Manifolds}

\author{Y{\i}ld{\i}ray Ozan}
\address{Department of Mathematics, Middle East Technical University,
\newline
06531 Ankara, Turkey} \email{ozan@metu.edu.tr}

\date{\today}
\thanks{The author is partially supported by the Turkish Academy of
Sciences (TUBA-GEBIP-2004-17).} \subjclass{Primary 57S05.
Secondary 53D35} \keywords{Trace homomorphism, Diffeomorphism
group.} \pagenumbering{arabic}

\begin{abstract}
Let $M$ be a closed connected smooth manifold and
$G=\textmd{Diff}_0(M)$ denote the connected component of the
diffeomorphism group of $M$ containing the identity. The natural
action of $G$ on $M$ induces the trace homomorphism on homology.
We show that the image of trace homomorphism is annihilated by the
subalgebra of the cohomology ring of $M$, generated by the
characteristic classes of $M$. Analogously, if $J$ is an almost
complex structure on $M$ and $G$ denotes the identity component of
the group of diffeomorphisms of $M$ preserving $J$ then the image
of the corresponding trace homomorphism is annihilated by
subalgebra generated by the Chern classes of $(M,J)$.
\end{abstract}

\maketitle

\section{Introduction and the results}
Let $G$ be any topological group acting on a topological space $X$
and $\it R$ any commutative ring. We define the trace
homomorphism,
$$H_k(G,\it R)\times H_l(X,\it R)
\stackrel{tr_*}{\rightarrow} H_{k+l}(X,\it R),$$ corresponding to
this action as follows: if $\phi:U\rightarrow G$ and
$\sigma:A\rightarrow X$ are cycles in $G$ and $X$ of degrees $k$
and $l$ representing classes $\upsilon$, $\alpha$, respectively,
let $tr_*(\upsilon,\alpha)$ be the class represented by the
homology cycle $(u,a)\mapsto \phi(u)(a)$, $(u,a)\in U\times A$. In
2003, it is proved in \cite{LM,LMP} that the trace homomorphism of
the Hamiltonian group of a closed symplectic manifold $(M,\omega)$
on the rational homology of $M$,
$$H_k(\textmd{Ham}(M,\omega),{\mathbb Q})\times H_l(M,{\mathbb Q})
\stackrel{tr_*}{\rightarrow}H_{k+l}(M,{\mathbb Q}),$$ is trivial,
for $k\geq 1$. Inspired by this result we prove the following
smooth analog:

\begin{Theorem}\label{thm-diff-trace}
Let $M$ be a closed connected smooth manifold and $\it R$ denote
the either field ${\mathbb Z}_2$ or ${\mathbb Q}$. Also let $P$
denote the subalgebra of the cohomology algebra $H^*(M,\it R)$,
generated by the Stiefel-Whitney classes $w_i(M)$, if $\it
R={\mathbb Z}_2$, and the subalgebra generated by the Pontryagin
classes $p_i(M)$ and the Euler class $e(M)$, if $\it R={\mathbb
Q}$. If $\textmd{Diff}_0(M)$ denotes the connected component of
the diffeomorphism group of $M$ containing the identity then the
image of the trace homomorphism
$$H_k(\textmd{Diff}_0(M),\it R)\times H_l(M,\it R)
\stackrel{tr_*}{\rightarrow} H_{k+l}(M,\it R)$$ is in the
annihilator of $P$, provided that $k\geq 1$.
\end{Theorem}

The proof of the above result yields immediately the following
almost complex analog:

\begin{Theorem}\label{thm-alcx-trace}
Assume that $M$ is a closed connected smooth manifold and $J$ is
an almost complex structure on $M$. Let $P$ denote the subalgebra
of the cohomology algebra $H^*(M,{\mathbb Q})$, generated by the
Chern classes $c_i(M)$. If $\textmd{Diff}_0(M,J)$ denotes the
identity component of the group of diffeomorphisms of $M$
preserving $J$, then the image of the trace homomorphism
$$H_k(\textmd{Diff}_0(M,J),{\mathbb Q})\times H_l(M,{\mathbb Q})
\stackrel{tr_*}{\rightarrow} H_{k+l}(M,{\mathbb Q})$$ is in the
annihilator of $P$, provided that $k\geq 1$.
\end{Theorem}

\subsection{Trace homomorphism on cohomology} For $\it R={\mathbb
Z}_2$ or ${\mathbb Q}$ we have $H^p(M,\it
R)=\textmd{Hom}(H_p(M,\it R),\it R)$ and using this duality we may
define trace homomorphism in cohomology: Let $u\in
H_k(\textmd{Diff}_0(M),\it R)$ and define
$$tr_u^*:H^p(M,\it R)\rightarrow H^{p-k}(M,\it R)$$ by the formula
$a\mapsto tr_u^*(a)$, $a\in H^p(M,\it R)$, where
$$tr_u^*(a):H_{p-k}(M,\it R)\rightarrow {\it R}, \
tr_u^*(a)(\alpha)=a(tr_*(u,\alpha)), \ \alpha \in H_{p-k}(M,\it
R).$$ Hence, the conclusions of Theorem~\ref{thm-diff-trace} and
of Theorem~\ref{thm-alcx-trace} can be written as $tr_u^*(P)=0$,
for all $u\in H_k(\textmd{Diff}_0(M),\it R)$, $k\geq 1$.

Suppose that $u\in H_k(\textmd{Diff}_0(M),\it R)$ is a spherical
class. Using any cycle representing $u$ we can build a fiber
bundle $M\rightarrow E\rightarrow S^{k+1}$, such that the
connecting homomorphism in the Wang sequence corresponding to this
bundle is nothing but the trace homomorphism:
$$\rightarrow H^{p-1}(E,\it R)\rightarrow H^{p-1}(M,\it R)
\stackrel{tr_u^*}{\rightarrow} H^{p-k-1}(M,\it R)\rightarrow
H^p(E,\it R)\rightarrow$$

It is well known that the connecting homomorphism in the Wang
sequence is a derivation of degree $k$ (\cite{Wh}). In other
words, for any $x,y \in H^*(M,\it R)$,
$$tr_u^*(xy)=tr_u^*(x) \ y+(-1)^{k \deg (x)}x \ tr_u^*(y).$$

On the other hand, for general $u$, since $\textmd{Diff}_0(M)$ is
an $H$-space any rational homology class $u$ is a product of
spherical classes (cf. see Section 5 of \cite{LM}) and therefore
$tr_u$ is the composition of the trace homomorphisms corresponding
to the spherical factors of $u$. Hence, we obtain the following
result:

\begin{Proposition}\label{prop-trace}
Let $u\in H_k(\textmd{Diff}_0(M),{\mathbb Q})$, $k>0$. For
cohomology classes $x,y \in H^*(M,{\mathbb Q})$ such that $y\in P$
(hence $tr_u^*(y)=0$) we have $tr_u^*(xy)=tr_u^*(x)y$. Moreover,
if $\deg (x)< k$ then $tr_u^*(xy)=0$.
\end{Proposition}

The above proposition yields the following corollary:

\begin{Corollary}\label{cor-module}
The natural map
$$tr^*:H_k(\textmd{Diff}_0(M),{\mathbb Q})\rightarrow \hom_P(H^*(M,{\mathbb
Q}),H^{*-k}(M,{\mathbb Q}))$$ is a homomorphism, where we regard
$H^*(M,{\mathbb Q})$ as a right module over its subalgebra $P$ and
$\hom_P(H^*(M,{\mathbb Q}),H^{*-k}(M,{\mathbb Q})$ denotes the
group of $P$-modulo homomorphisms.
\end{Corollary}

\begin{Example}\label{exm-Kahler}
Let $u$ be as in the complex analog of the above proposition,
where $P$ is generated by the Chern classes of the almost complex
manifold $(M,J)$ and $u$ belongs to $H_k(\textmd{Diff}_0(M,J),\it
R)$. Assume that $(M,\omega)$ is a monotone closed symplectic
manifold of dimension $2n$. So $[\omega]$ is a multiple of
$c_1(M)$ and hence it is in $P$. Assume further that $M$ has the
Hard Lefschetz Property, i.e., $$\cup \ [\omega]^r
:H^{n-r}(M,{\mathbb C})\rightarrow H^{n+r}(M,{\mathbb C})$$ is an
isomorphism for any $r\geq 0$. So, if $b\in H^{n+r}(M,{\mathbb
C})$ then $b=a \ [\omega]^r$ for some $a\in H^{n-r}(M,{\mathbb
C})$ and hence $$tr_u^*(b)=tr_u^*(a \ [\omega]^r)=tr_u^*(a) \
[\omega]^r.$$ In particular, $tr_u^*([\omega]^r)=0$. It follows
that, if $k>n$ then $tr_u^*=0$.
\end{Example}

\section{Proof of the Theorem}
To prove the above results we need to recall the definition and
some basic properties of equivariant bundles: Let $G$ be any Lie
group and $F\rightarrow E\stackrel{\pi}{\rightarrow} B$ a fiber
bundle. If $G$ acts on both $E$ and $B$ such that the projection
map $\pi$ is $G$-equivariant; i.e., $\pi (v\cdot g) = \pi (v)
\cdot g$, for all $g\in G$ and $v \in E$, we say that the bundle
is $G$-equivariant. Note that if $X$ is also a $G$-space and
$f:X\rightarrow B$ is a $G$-equivariant map then the pullback
bundle has an induced $G$-equivariant structure.

\begin{Example}\label{exm-Equiv-bundle}
{\bf i)} Let $F\rightarrow E\stackrel{\pi}{\rightarrow} B$ be a
G-equivariant fiber bundle, where the action of $G$ on $B$, and
hence on $E$, is free.  Taking quotients of both the total space
and the base by $G$, we get another fiber bundle $F\rightarrow
E/G\stackrel{\tilde{\pi}}{\rightarrow} B/G$, whose pullback via
the quotient map $p:B\rightarrow B/G$ is isomorphic to the bundle
$F\rightarrow E\stackrel{\pi}{\rightarrow} B$.

{\bf ii)} Let $M$ be a smooth manifold.  Since any diffeomorphism
of $M$, $\phi:M\rightarrow M$, extends to the tangent bundle
$\phi_*:T_*M\rightarrow T_*M$, we see that the tangent bundle is
$\textmd{Diff}(M)$-equivariant, where $\textmd{Diff}(M)$ is the
group of all diffeomorphisms of $M$.
\end{Example}

\begin{proof}[Proof of Theorem~\ref{thm-diff-trace}]
Let $G$ denote the group $\textmd{Diff}_0(M)$, the group of
diffeomorphism of $M$ isotopic to the identity, and
$\sigma:A\rightarrow M$ be a cycle in $M$ representing any given
class $\alpha$ of degree $l$. Since the base field is either
${\mathbb Z}_2$ or ${\mathbb Q}$ we may assume that $A$ is a
closed smooth manifold and $\sigma:A\rightarrow M$ is a smooth
map. Consider the trace map
$$tr:A\times G \rightarrow M, \hspace{0.2cm} (a,g)\mapsto
\sigma(a)\cdot g,  \ \mbox{for all} \ (a,g)\in A\times G.$$  To
prove the theorem it suffices to show that $tr^*(v)=0$, for any
$v\in P$ of degree $l+k$.

Note that $G$ acts on $A\times G$ by right multiplication on the
second factor, which makes the trace map $G$-equivariant. By the
above example the tangent bundle $T_*M\rightarrow M$ is
$G$-equivariant and hence the pullback bundle
$tr^*(T_*M)\rightarrow A\times G$ is $G$-equivariant. Since the
$G$-action on the base space $A\times G$ is free this bundle is
induced from the quotient bundle $tr^*(T_*M)/G\rightarrow (A\times
G)/G$, which is isomorphic to $\sigma^*(T_*M)\rightarrow A$. In
particular, by the naturality of characteristic classes $tr^*(v)$
is the pullback of a class in $H^{k+l}(A,\it R)$.  However,
$H^{k+l}(A,{\it R})=0$, because $A$ is $l$-dimensional and $k\geq
1$. Hence, $tr^*(v)=0$ and the proof finishes.
\end{proof}

\begin{Remark} {\bf i)} Note that the above proof works also for
Theorem~\ref{thm-alcx-trace}. Indeed more is true: Let $G$ be a
subgroup of $\textmd{Diff}_0(M)$ and $E\rightarrow M$ be a
$G$-equivariant real or complex vector bundle. Then the analogous
result to Theorem~\ref{thm-diff-trace} holds for $G$ and the
subalgebra $P$ of the cohomology algebra of $M$, generated by the
characteristic classes of $E$.

{\bf ii)} Another extension of the main theorem, suggested by
Dieter Kotschick, to foliations is as follows: Assume that the
smooth manifold $M$ is foliated and let $G$ be the subgroup of
$\textmd{Diff}_0(M)$ preserving the foliation. Then $G$ acts as
isomorphisms of both the tangent bundle $\tau$ and the normal
bundle $\eta$ of the foliation, where we have
$T_*M=\tau\oplus\eta$. Then for this $G$ we can replace the
subalgebra $P$ in Theorem~\ref{thm-diff-trace} with the subalgebra
generated by the characteristic classes of the two summands of the
tangent bundle to $M$. Note that this subalgebra is clearly bigger
than the subalgebra generated by the characteristic classes of
$T_*M$ only. Of course, this is no surprise since $G$ is generally
much smaller than $\textmd{Diff}_0(M)$.
\end{Remark}

\subsection*{Acknowledgment} The author would like to thank Dusa McDuff
for pointing out some errors in the earlier version of the article
and to Dieter Kotschich for his remarks and suggestions.

\providecommand{\bysame}{\leavevmode\hboxto3em{\hrulefill}\thinspace}


\begin{thebibliography}{1}

\bibitem{LM} {\sc F. Lalonde, D. McDuff}. Symplectic structures on fiber
bundles. {\it Topology} {\bf 42} (2003), 309-347.

\bibitem{LMP} {\sc F. Lalonde, D. McDuff, L. Polterovich}. Topological rigidity of
Hamiltonian loops and Quantum homology. {\it Invent. Math.} {\bf
135} (1999), 369-385.

\bibitem{Wh} {\sc G. W. Whitehead}. {\it Elements of Homotpy Theory}. 3rd Edition,
 (Springer-Verlag, 1995) pp. 319-320.

\end{thebibliography}
\end{document}